\documentclass{elsarticle}

\usepackage{amssymb}
\usepackage{graphicx}

\newtheorem{theorem}{Theorem}
\newtheorem{proposition}{Proposition}
\newtheorem{lemma}{Lemma}

\newtheorem{example}{Example}

\newtheorem{property}{Property}
\newtheorem{remark}{Remark}

\newcommand{\R}{\mathbb{R}}

\newcommand{\Q}{\mathbb{Q}}

\journal{Journal of Number Theory}

\begin{document}

\begin{frontmatter}

\title{Polygons in billiard orbits}
\author{Henk Don}
\ead{henkdon@gmail.com}
\address{TU Delft, EWI (DIAM), Section Probability Theory, Mekelweg 4, 2628 CD Delft, the Netherlands. Phone: +31 152782546}

\begin{abstract}
We study the geometry of billiard orbits on rectangular billiards. A truncated billiard orbit induces a partition of the rectangle into polygons. We prove that thirteen is a sharp upper bound for the number of different areas of these polygons.   
\end{abstract}

\begin{keyword}
billiard orbit \sep geometry of partitions
\MSC[2010] 11B75
\end{keyword}

\end{frontmatter}

\section{Introduction}

Let a billiard ball be shot from a corner of a rectangular billiard. Consider the ball as a point, and truncate the orbit somewhere at the boundary. The truncated orbit of the ball generates a partition of the rectangular billiard into polygons, similar to Figure \ref{billiard}. Many of these triangles and quadrangles seem to have the same shape and size. In this paper we will show that (for a fixed shooting angle and stopping point) the number of different areas is at most thirteen. This universal upper bound is the sharpest possible. We also consider rational shooting angles and irrational shooting angles for which the thirteen is never reached.

\begin{figure}[!h]
\begin{center}
\includegraphics*[width =6.0cm]{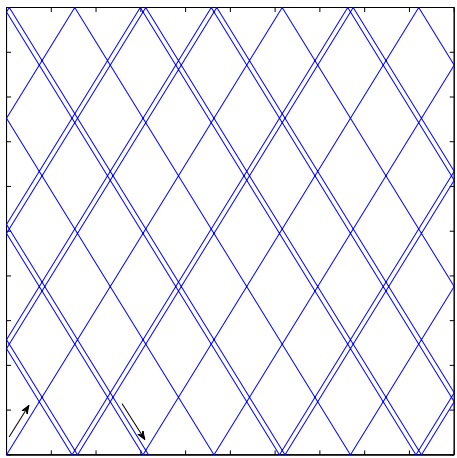}
\caption{ \footnotesize Truncated orbit of a billiard ball. The arrows indicate start and end of the orbit.}
\end{center}
\end{figure}\label{billiard}

\section{Rotations}

\begin{figure}[!h]
\begin{tabular}{r|l}
\includegraphics*[width =5.8cm]{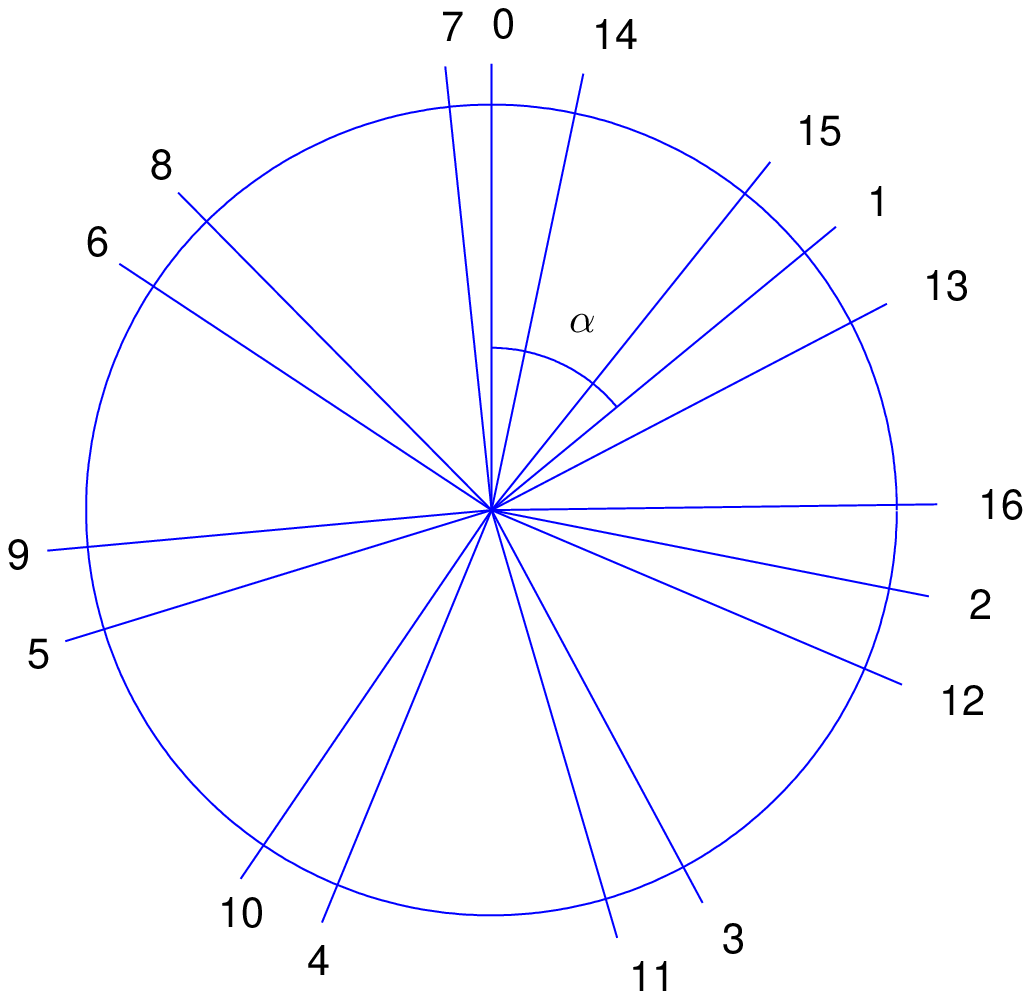}
&
\includegraphics*[width =5.8cm]{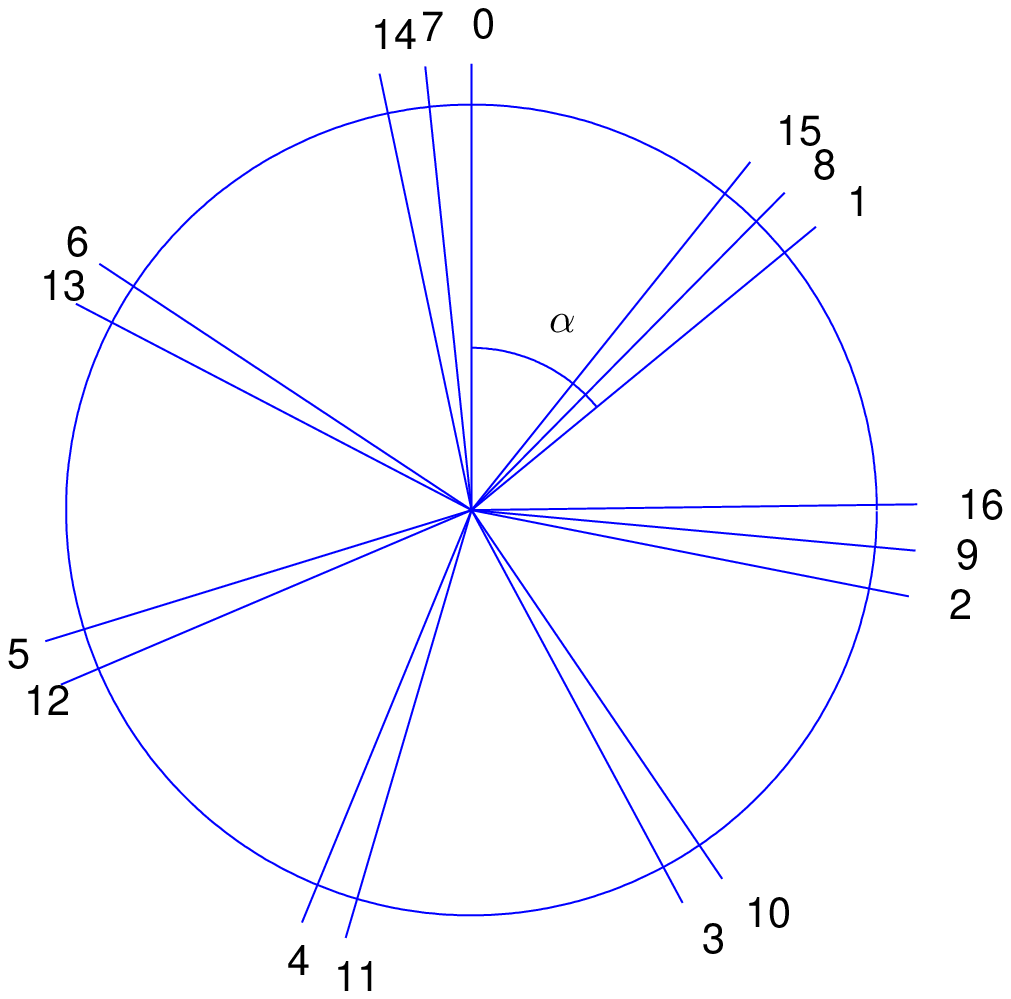}
\end{tabular}\caption{Left figure, the Three Gap Theorem: Cutting
a pie $n$ times where each next cut is obtained by shifting the previous one over a fixed angle $\alpha$ gives
at most three different sizes of pieces of the pie. Right figure, the Four Gap Theorem: Now the first cut (at
$0$) works as a `reflecting boundary'. As soon as it is reached, we continue in the opposite direction. In this
case we have after $n$ cuts at most four different sizes. \newline For this picture we used $\alpha =
0.1405*2\pi$ and $n=17$.}\label{3GT4GT}
\end{figure}

The results in this paper are closely related to the Three Gap Theorem (see e.g. \cite{Sos}, \cite{Rav}) and the Four Gap
Theorem (see \cite{Don}). The statements of these two theorems are best illustrated by a picture; see Figure
\ref{3GT4GT}.

The Three Gap Theorem is naturally associated to the concept of \emph{rotations}. First we recall the theorem and
then we discuss rotations on intervals. For $x\in\mathbb{R}$, let $\left\{x\right\} = x -\lfloor x\rfloor$ denote its fractional part.
\begin{theorem} (The Three Gap Theorem)
Let $n\in\mathbb{N}$ and $\alpha \in (0,1)$. The numbers
\begin{equation}\label{rotation_unit}
0,\left\{\alpha\right\},\left\{2\alpha\right\},\left\{3\alpha\right\},\ldots,\left\{n\alpha\right\}
\end{equation}
induce a partition of the interval $[0,1]$ in subintervals which can have at most three different lenghts. If there are three lengths, then the largest is the sum of the other two.
\end{theorem}
Letting $T_\alpha(x) = \left\{x+\alpha\right\}$ for $x\in [0,1]$, the numbers (\ref{rotation_unit}) transform into
\begin{equation}\label{rotation_orbit}
0,T_\alpha(0),T^2_\alpha(0),\ldots,T^n_\alpha(0).
\end{equation}
If we consider $x$ as a point on the circle of unit circumference, then $T_\alpha(x)$ is obtained by rotating $x$ over a distance $\alpha$. This gives a more dynamical view of the partition of $[0,1]$: the partition is induced by a truncated orbit of the rotation map $T_\alpha$. These observations lead to the following generalization of the Three Gap Theorem:
\begin{property}\label{prop_rotation}
Let $n_1,n_2\in\mathbb{N}$, $\alpha \in (0,1)$ and $a,b\in\mathbb{R}$. The $n_1+n_2+1$ numbers
\begin{equation}\label{rotation_general}
aT_\alpha^{-n_1}(0)+b, \ldots,aT_\alpha^{-1}(0)+b,\,b\,,aT_\alpha(0)+b,\ldots,aT_\alpha^{n_2}(0)+b
\end{equation}
induce a partition of $[b,b+a]$ in subintervals having at most three different lenghts.
\end{property}
This can easily be obtained by taking $n = n_1+n_2$ in (\ref{rotation_orbit}), rotating over an appropriate angle and applying the linear map $a\cdot+b$ to the orbit. Actually, a special case of this property already appeared as a theorem in \cite{Don}. However, there a
complicated proof was given to obtain this result. Vilmos Komornik came up with the idea to place the numbers on
the circle, thus obtaining a much simplified and more natural argument \cite{Komornik}. In the sequel we will refer to (\ref{rotation_general}) as a \emph{rotation orbit on} $[b,b+a]$.\\
There is a slightly stronger property we will need in Remark \ref{remark} (see e.g. \cite{Don} and \cite{Rav}):
\begin{property}\label{prop_strong}
Take a truncated orbit of a rotation on an interval. Suppose the orbit consists of $n$ numbers. Create another orbit from this by removing the last number. The two partitions induced by these orbits give two sets of lengths. The union of these two sets contains at most three different lengths.
\end{property}

\section{Billiards and the Four Gap Theorem}

\begin{figure}[!h]
\begin{center}
\includegraphics*[width =9cm]{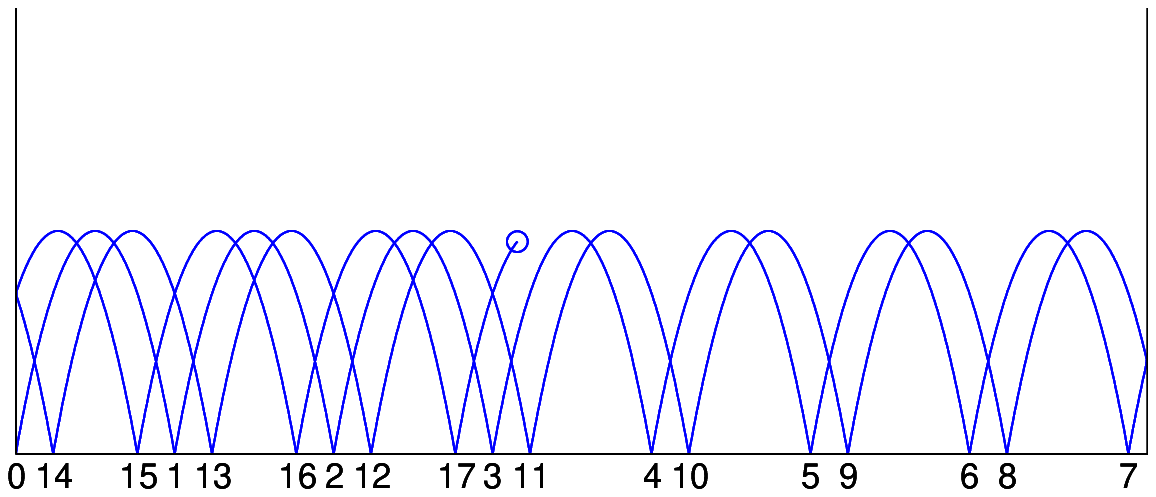}
\caption{\footnotesize A ball bouncing on an interval between two walls. The Four Gap Theorem makes a statement
about the subset of the interval consisting of the landing points of the ball. We used $0.1405$ times the length
of the interval as bouncing distance.}\label{1D_bounce}
\end{center}
\end{figure}

The billiard in Figure \ref{billiard} can be seen as a generalization to two dimensions of the pie-cutting
process of the Four Gap Theorem, as illustrated in Figure \ref{3GT4GT}. This statement deserves some
explanation. Figure \ref{1D_bounce}, a picture in some sense equivalent to the right panel of Figure
\ref{3GT4GT}, gives a description of the Four Gap Theorem in terms of a ball bouncing on the unit interval.

\begin{figure}[!h]
\begin{center}
\includegraphics*[width =9cm]{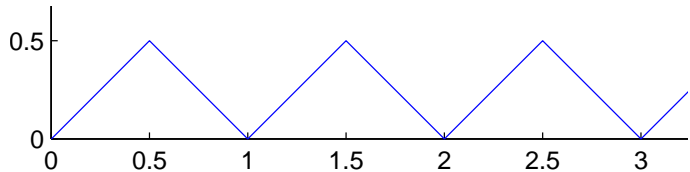}
\caption{\footnotesize Plot of the folding map $f_1(x) = ||x||$.}\label{folding}
\end{center}
\end{figure}

This figure shows the movement of a ball bouncing between two walls, where we assume that the ball is a point
and that there is no loss of energy. The landing points of the ball build a sequence in the interval. The first
$n$ numbers in this sequence ($0$ included) define a splitting of the interval in $n$ subintervals. The main
statement of the Four Gap Theorem is that these subintervals can have at most four different lengths. In Figure
\ref{billiard} we now have a subset of a square, consisting of those points where the billiard ball
appears. This observation gives already some reason to consider the billiard as a $2$-dimensional generalization
of the pie of the Four Gap Theorem. However, we can also argue this point
of view in a more mathematical way.\\
Let $||x||$ denote the distance from $x$ to the nearest integer. For $\alpha\in\R\setminus\Q$, let
$$
S_\alpha := (||k\alpha||)_{k=0}^\infty.
$$
Obviously this is a sequence in $[0,\frac{1}{2}]$. Moreover, it is exactly the sequence of landing points of a
ball bouncing between $0$ and $\frac{1}{2}$ with horizontal bouncing distance $\alpha$. The sequence $S_\alpha$
is obtained by `folding' the sequence of integer multiples of $\alpha$ into the interval $[0,\frac{1}{2}]$. What
we mean by this folding is illustrated in Figure \ref{folding}, where we plot the function
$f_1:[0,\infty)\rightarrow[0,\textstyle\frac{1}{2}]$
$$
f_1(x) := ||x||,
$$
and illustrate how $[0,\infty)$ is mapped to $[0,\frac{1}{2}]$ by $f_1$.

Now we concentrate on the billiard: the orbit of the billiard ball is obtained by `folding' a halfline into a
rectangle. Since the shooting angle is arbitrary between $0$ and $\pi/2$, we may without loss of generality
assume that instead of a rectangle the billiard is a square and equal to $[0,\frac{1}{2}]^2$. The `folding' map
corresponding to this billiard is given by a two-variable function $f_2:[0,\infty)^2
\rightarrow[0,\textstyle\frac{1}{2}]^2$:
$$
f_2(x,y) = (||x||,||y||).
$$
As we see, $f_2(x,y) = (f_1(x),f_1(y))$, which is why the billiard can be viewed as being a generalization of
the setting of the Four Gap Theorem to two dimensions. The folding map $f_2$ applied to a line creates a
billiard orbit. Let $\alpha>0$, then
$$
B_{[0,M]}^\alpha := \left\{(||x||,||\alpha x||):x\in [0,M]\right\}
$$
describes a truncated billiard orbit that has initial slope $\alpha$ (the slope alternates between $\alpha$ and
$-\alpha$). Let $\mathcal{A}_{[0,M]}^\alpha$ and $\mathcal{S}_{[0,M]}^\alpha$ denote the number of different
areas respectively different shapes in the partition of $[0,\frac{1}{2}]^2$ induced by $B_{[0,M]}^\alpha$. Two
shapes are different if one can not be obtained from the other by translating, rotating and reflecting.
For orbits truncated in a boundary point, we will prove the following theorem:

\begin{theorem}\label{theorem_13areas}
Let $\alpha>0$ and choose $M>0$ such that $(||M||,||\alpha M||)\in
[0,\frac{1}{2}]^2\setminus(0,\frac{1}{2})^2$. Then the billiard orbit $B^\alpha_{[0,M]}$ induces a partition of
$[0,\frac{1}{2}]^2$ in polygons for which
$$
\mathcal{A}_{[0,M]}^\alpha \leq 13\quad \textrm{and}\quad \mathcal{S}_{[0,M]}^\alpha \leq 16.
$$
These upper bounds are the best possible.
\end{theorem}

In this theorem the billiard is square, but the result for rectangular billiards easily follows since the square
can be scaled to any rectangle without changing the ratios between the shapes. From now on, we will assume that $M$ satisfies the condition in the theorem.

\section{Orbit construction}

\begin{figure}[!h]
\begin{tabular}{r|c|l}
\includegraphics*[width =3.7cm]{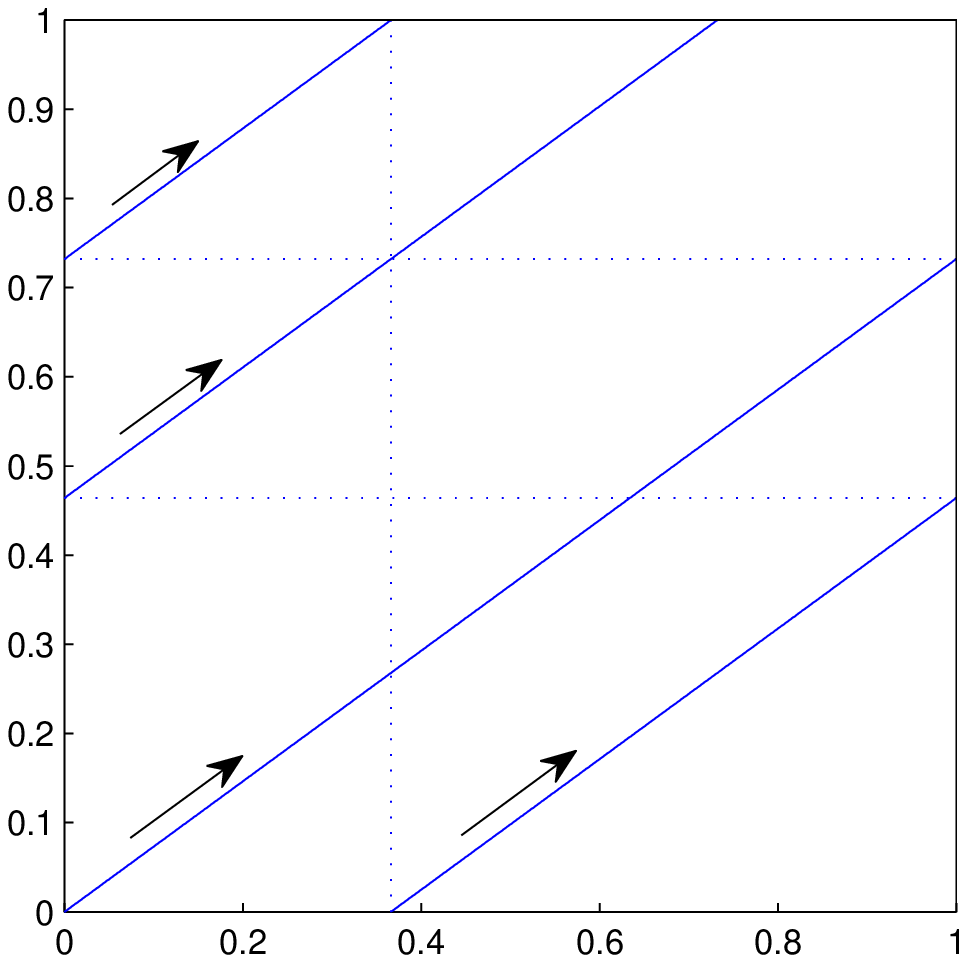}
&
\includegraphics*[width =3.7cm]{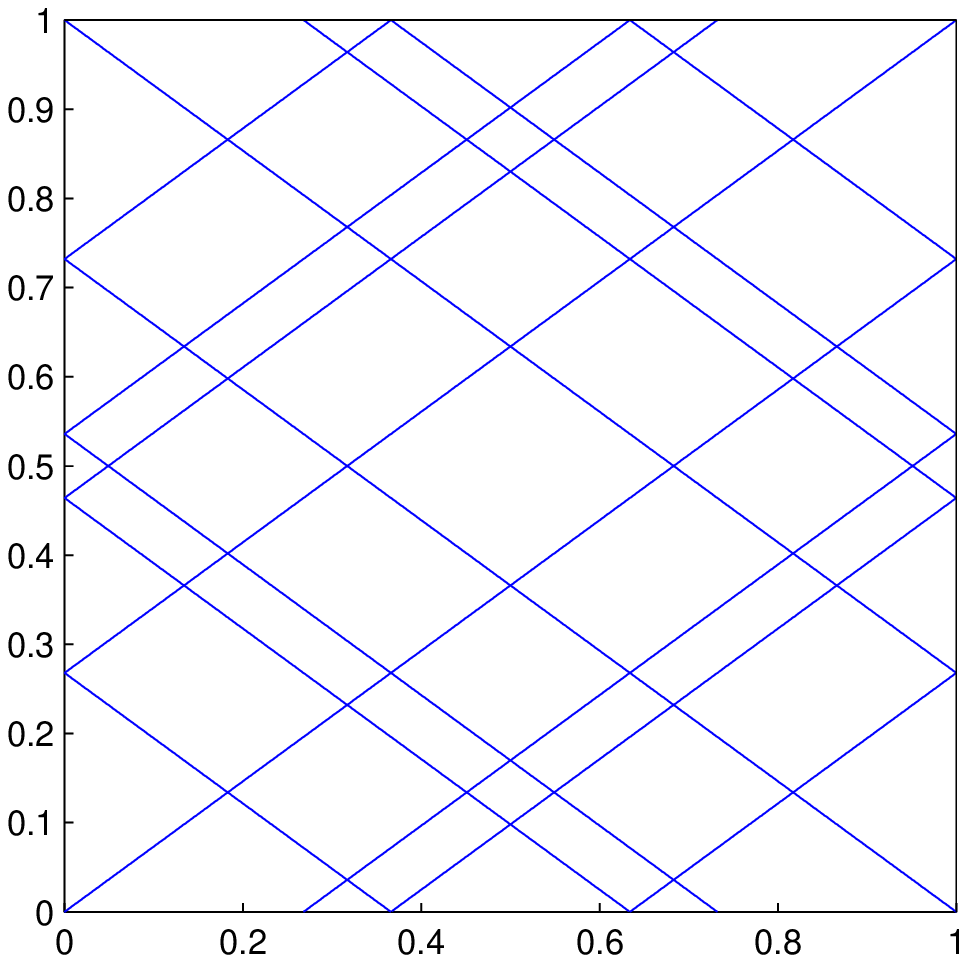}
&
\includegraphics*[width =3.7cm]{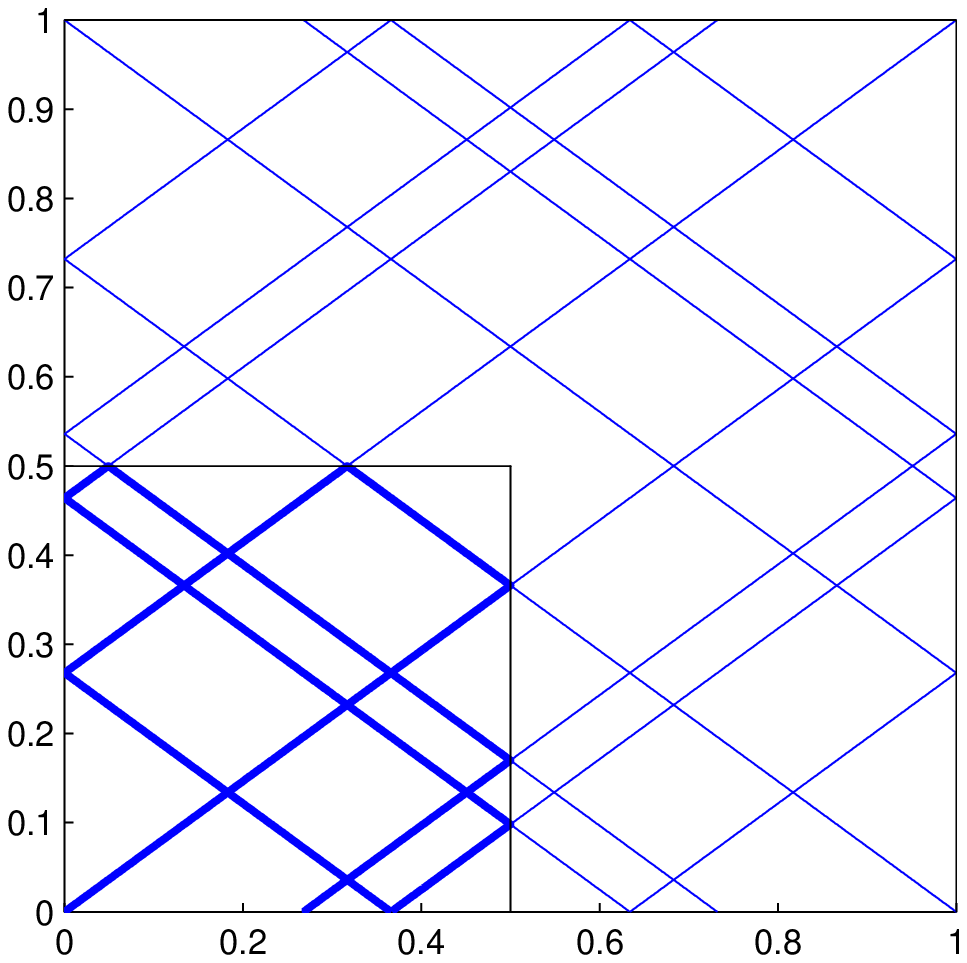}
\end{tabular}\caption{
Construction of a billiard orbit in three steps. Here $N=4$ and $\alpha = \sqrt{3}-1$.}\label{construction}
\end{figure}

We already have an explicit expression for the billiard orbit $B_{[0,M]}^\alpha$, but we will need a more
tractable description. Therefore, in this section we present a rough intuitive outline of the way one can think
of the geometry and construction of the billiard. The corresponding lemmata and their proofs are given in Section
\ref{section_proofs}. Consider the unit square and draw a line starting from the lower left corner with slope
$\alpha$. The boundaries are now considered to be connected as in a torus, so when we reach it, the line
continues at the opposite boundary. Equivalently, if one of the coordinates is about to exceed $1$, we subtract
$1$. But this is exactly taking fractional parts in both coordinates. Therefore, after we have traversed the
unit square $N$ times, we have a set which can be expressed as
$$
\left\{(\left\{x\right\},\left\{\alpha x\right\}):0\leq x < M\right\},
$$
for some $M\in\mathbb{R}$. A plot of such a set is shown in the left panel of Figure \ref{construction}.

Now do the same starting from the other corners, traversing the square $N$ times with a line either with slope
$\alpha$ or $-\alpha$. Explicit expressions for these four sets (one for each corner) are given in Lemma \ref{lemma_Aexpressions}.
For an illustration, see the middle plot in Figure \ref{construction}.\\
The key observation now is that intersection of all $4N$ lines with $[0,\frac{1}{2}]^2$ gives exactly a truncated billiard orbit with slope $\alpha$, as is proved in Lemma \ref{lemma_orbit}. This fact is illustrated in the right plot in Figure \ref{construction}. Obviously not all $4N$ lines actually contribute to the billiard orbit. However, there is a good reason to consider them all: the intercepts of the $2N$ lines with positive slope form a truncated orbit of a rotation on the interval $[-\alpha,1]$, see Lemma \ref{lemma_rotation}. For the lines with negative slope a similar result holds. Having collected these insights, a simple counting argument suffices to obtain the upper bounds claimed in Theorem \ref{theorem_13areas}:\\
\\
\textbf{Proof of Theorem \ref{theorem_13areas}} Lemma \ref{lemma_orbit} writes the billiard orbit as an intersection of the square $[0,\frac{1}{2}]^2$ with a set of lines. Let us concentrate on the lines with positive slope. By Lemma \ref{lemma_rotation} the intercepts of these lines form a rotation orbit on the interval $[-\alpha,1]$. So by Property \ref{prop_rotation} they induce a partition of this interval in subintervals of \textit{at most} three different lengths. Denote the set of these lengths by $D :=
\left\{d_1,\ldots,d_n\right\}$, where $n\leq 3$. For the lines with negative slope, the intercepts are the numbers $1-y_k$, $-N\leq
k\leq N$. They induce a partition of $[0,1+\alpha]$ in subintervals having lengths in the same set $D$. It now follows that vertical distances between adjacent parallel lines are in the set $D$.

\begin{figure}[!h]
\begin{center}
\includegraphics*[width =6cm]{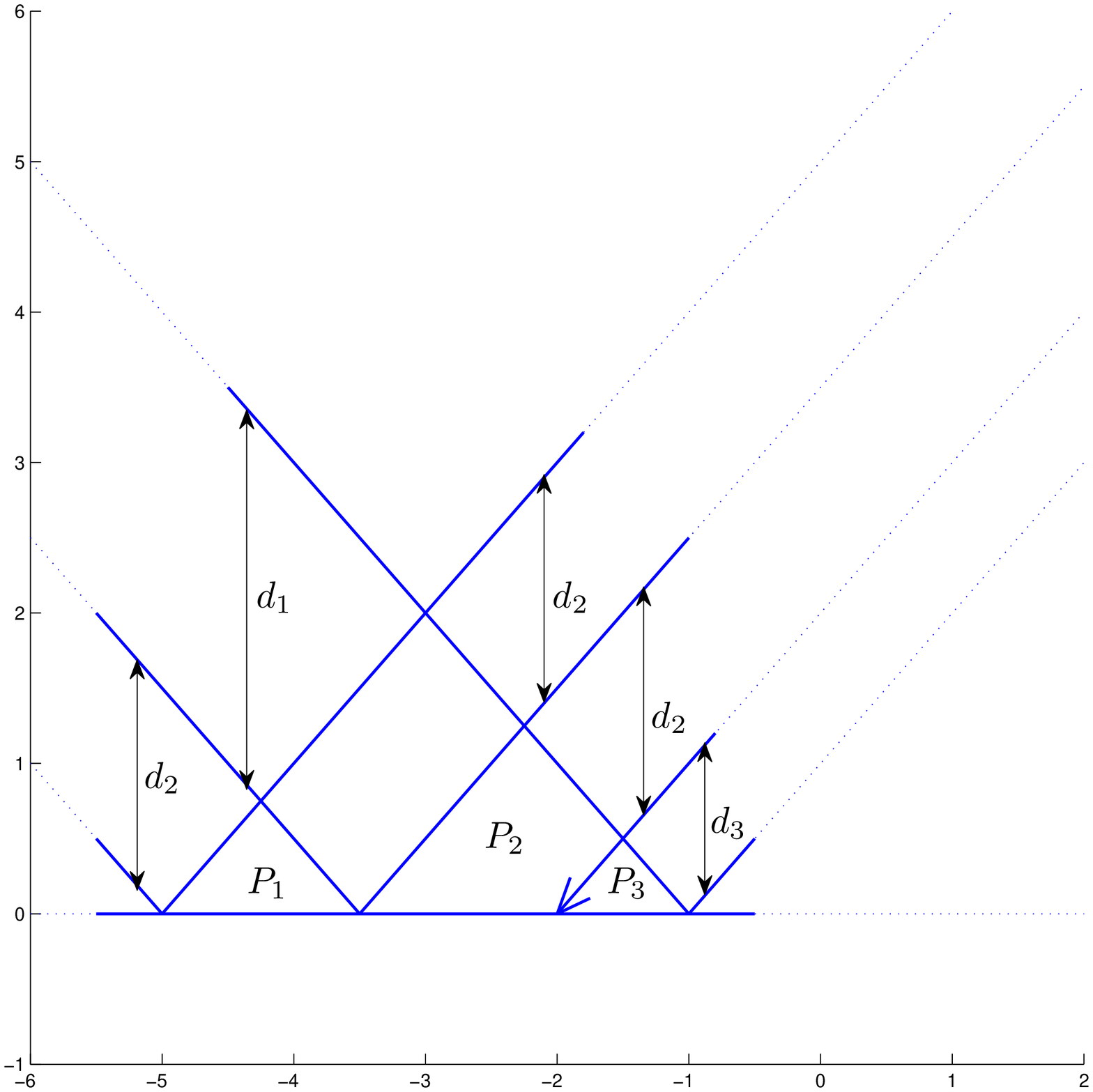}
\caption{\footnotesize Local situation at the boundary where the orbit ends. Polygons of type $2$ are
triangular if the endpoint of the orbit is not one of the corners of the polygon (as is the case with $P_1$).
There are only two shapes for which the endpoint of the orbit is one of the corners. One of them is still
triangular (in this example $P_3$), the other is irregular ($P_2$).}\label{pic_extreme}
\end{center}
\end{figure}

We will distinguish between three types of polygons: those that have no side which is part of the boundary of $[0,\frac{1}{2}]^2$ (type $1$), those that have exactly one such a side (type $2$) and those that have two or more (type $3$).\\
The polygons of type $1$ must be parallelograms. The
area of such a parallelogram is given by $d_id_j/2\alpha$ for some $d_i,d_j\in D$, and consequently they can have at most six different areas.\\
A polygon of type $2$ that is triangular must be half of a rhombus of which the vertical diagonal has length $d\in
D$, and therefore its area is $d^2/4a$. There is at most one non-triangular type $2$ polygon, as is explained in Figure \ref{pic_extreme}. So polygons of type $2$ can have at most four different areas.\\
Polygons of type $3$ must be in one of the corners of $[0,\frac{1}{2}]^2$, but not in $(0,0)$ since the orbit starts there. So this gives at most three more areas.\\
Putting everything together, it turns out that the number of different areas is bounded by thirteen.\\
For the number of shapes a similar counting argument holds. The number of parallelogram shapes is again six, since reflections do not count. The triangles that are half of a rhombus can have at most six different shapes, since there are three types of rhombi which can be cut either horizontally or vertically. The rest of the argument doesn't change, so there are at most three more different shapes than different areas, which establishes the upper
bound of at most sixteen different shapes.\\
The sharpness of these bounds follows from Example \ref{ex_13flakes} in section \ref{example}. \hfill $\Box$

\begin{remark}\label{remark}\rm{As the careful reader may have noted, the construction of the billiard orbit always gives a truncation on the left boundary or on the lower boundary of the square. So strictly speaking, Theorem \ref{theorem_13areas} is not proved in full generality yet. Suppose we have an orbit truncated at the upper or right boundary. By removing the last linear part or adding the next linear part, we can transform this orbit into an orbit truncated at the left or lower boundary. This means that in the proof above, the rotation orbit on the interval $[-\alpha,1]$ contains one element more or one less than the rotation orbit on $[0,1+\alpha]$. Now Property \ref{prop_strong} tells us that vertical distances between adjacent parallel lines can still have at most three different values, completing the proof.}
\end{remark}

\section{Rational angles and a golden exception}

\begin{figure}[!h]
\begin{center}
\includegraphics*[width =6cm]{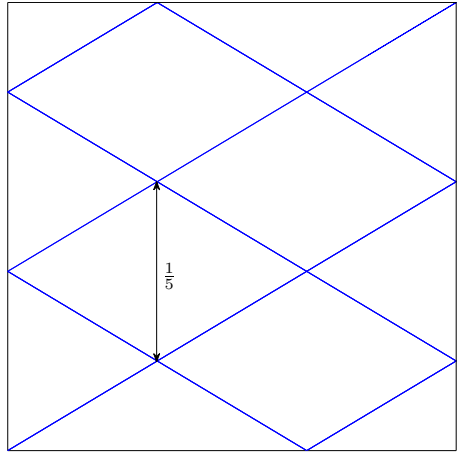}
\caption{\footnotesize The periodic orbit for $\alpha = 3/5$. There are three different areas: the rhombi have area $1/(2\cdot 3\cdot 5) = 1/30$. The triangles have area $1/60$ or $1/120$.}\label{pic_rational}
\end{center}
\end{figure}

Theorem \ref{theorem_13areas} gives an upper bound for the number of different areas of shapes on the billiard
table. Some natural questions remain. For example, what happens if $\alpha$ is rational? Can we prove sharper upper bounds under suitable conditions? In this section we
explore these properties.\\
Obviously, taking $\alpha$ rational gives a special case. The first thing to note is that the orbit will be periodic: if $\alpha = p/q$, then for $x\in\mathbb{R}$
$$
(||x+q||,||\alpha(x+q)||) = (||x||,||\alpha x||).
$$
A bit less trivial is the following result.
\begin{proposition}
The best upper bound for
$\mathcal{A}_{[0,M]}^\alpha$ with $\alpha\in\mathbb{Q}$ is $13$, but  for all $\alpha\in\mathbb{Q}$ there is an $M_0$ such that $1\leq
 \mathcal{A}_{[0,M]}^\alpha \leq 3$ for $M\geq M_0$. These bounds are sharp.
\end{proposition}
\textbf{Proof} Note that the areas of the polygons continuously depend on $\alpha$. So if we have an
$\tilde\alpha$ and $M$ such that $\mathcal{A}_{[0,M]}^{\tilde\alpha} = 13$, then we can find $\varepsilon > 0$
such that the upper bound of thirteen is reached for all $\alpha\in
(\tilde\alpha-\varepsilon,\tilde\alpha+\varepsilon)$. Since this interval contains rationals, we see that
rationality is not sufficient for a sharper upper bound.\\
Since the orbit is periodic, the partition doesn't change anymore if $M$ is large enough. Taking $\alpha = 1$ shows that $1$ is a sharp lower bound for the limiting number of shapes. For the upper bound, suppose that $\alpha = p/q$. By Lemma \ref{lemma_rotation} the intercepts satisfy
$$
y_{p+q} = \Bigl(1+ \frac{p}{q}\Bigr)\left\{\frac{(p+q)p/q}{1+p/q}\right\}-\frac{p}{q} = -\frac{p}{q} = y_0.
$$
It follows that the numbers $y_k$ form a periodic rotation orbit on $[-\alpha,1]$ and therefore the set $D$ as defined in the proof of Theorem \ref{theorem_13areas} contains only one length if $M$ is large enough. If $p$ and $q$ are relative prime, then this length is $1/q$. Now a type $1$ polygon is a rhombus with area $1/2pq$. Since there is no endpoint of the orbit anymore, a type $2$ polygon is half of such a rhombus. Polygons in the corners are also triangular, because the orbit touches all sides of the square before becoming periodic. These triangles are quarters of the rhombus, thus having area $1/8pq$. This makes at most three different areas in total. To see that this upper bound is sharp, see Figure \ref{pic_rational}.\hfill $\Box$

Surprisingly, there exist irrational $\alpha$ for which the upper bound of thirteen different areas is never
reached:

\begin{proposition}
Let $\phi = (\sqrt{5}-1)/2$ denote the small golden mean. If $\alpha = \frac{1}{n+\phi}$ for some $n\in\mathbb{N}$, then $\mathcal{A}_{[0,M]}^\alpha\leq 12$.
\end{proposition}

\textbf{Proof} Consider the numbers $y_k$ that form a rotation orbit on $[-\alpha,1]$. The partition of $[-\alpha,1]$ induced by this orbit gives subintervals with lengths in a set $D$. This set $D$ changes if we extend the orbit (i.e. we increase $M$): some lengths will disappear and new lengths will be created. In \cite{Don} and \cite{Rav} it was shown that the largest length is always the first to disappear. A new length only pops up if there are only two lengths in $D$, and the new length is the difference of these two existing lengths. Together with the fact that $1-\phi = \phi^2$, this is the basis of our argument.\\
Let $\alpha = 1/(n+\phi)$. From the way points are added to the rotation orbit it is clear that we can choose $M$ such that $[-\alpha,1]$  will be partitioned in $n+1$ intervals of length $\alpha$ and an interval of length $1+\alpha-(n+1)\alpha = \phi\alpha$. This gives $D = \left\{\alpha,\phi\alpha\right\}$. Extending the orbit with one more point transforms $D$ into $\left\{\alpha,\phi\alpha,\phi^2\alpha\right\}$ and this is the first time that $D$ contains three lengths. Increasing $M$ further, $D$ will change into $\left\{\phi\alpha,\phi^2\alpha\right\}$ and then into $\left\{\phi\alpha,\phi^2\alpha,\phi^3\alpha\right\}$. An inductive argument suffices to show that the ratios between the lengths in $D$ are preserved.\\
Recall that the areas of the parallelograms are determined by a product of two lengths in $D$. By the above reasoning, if $D = \left\{d_1,d_2,d_3\right\}$, then $d_1d_3 = d_2^2$, which implies that the parallelograms can have at most five different areas. Consequently $\mathcal{A}_{[0,M]}^\alpha\leq 12$.  \hfill$\Box$

\section{Lemmata and their proofs}\label{section_proofs}

\begin{figure}[!h]
\begin{center}
\includegraphics*[width =9cm]{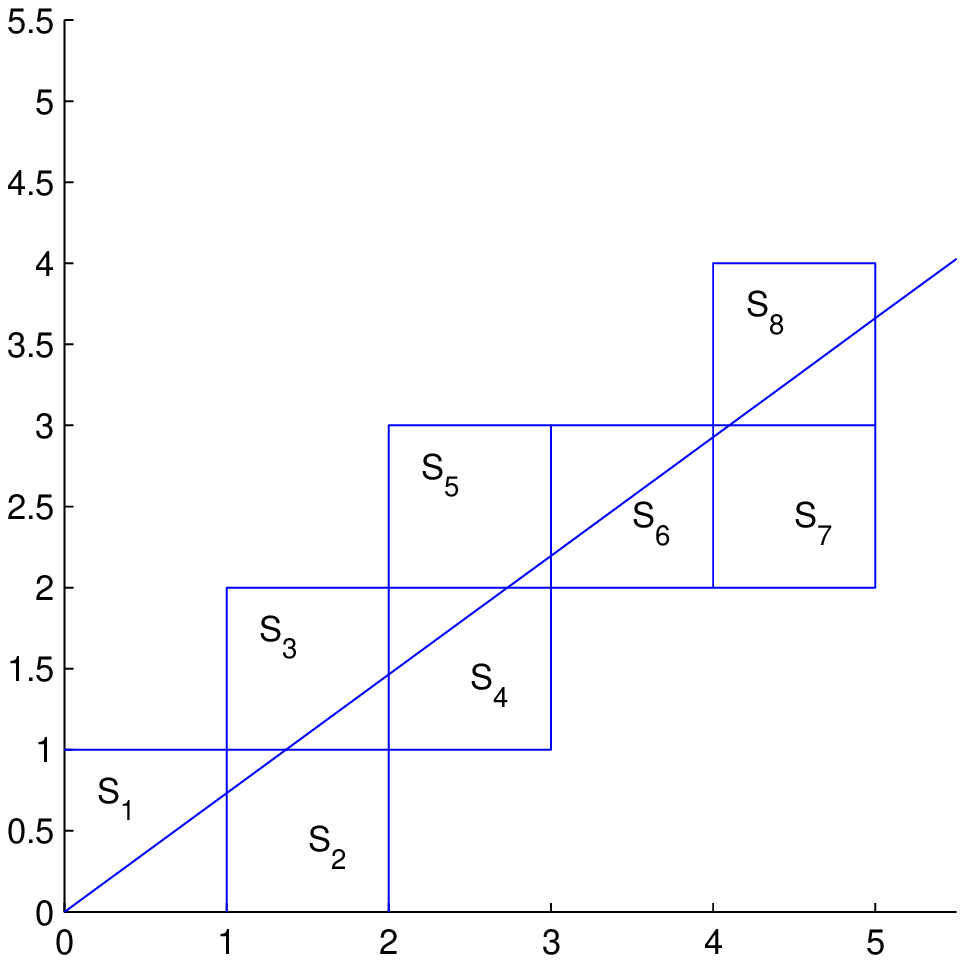}
\caption{\footnotesize Construction of the squares $S_k$. Here $\alpha = \sqrt{3}-1$, $N = 8$ and $M = 5$. The numbers $y_k$ are approximately given by $y_0 = 0$, $y_1 \approx 0.732$, $y_2\approx -0.268$, $y_3\approx 0.464$,\ldots Compare with Figure \ref{construction}, left plot.}\label{squares}
\end{center}
\end{figure}

Let $\alpha>0$ be an irrational number and consider the halfline $l(x) = \alpha x, x\geq 0$. Let $S_1 = [0,1)^2$ and define $S_2, S_3, S_4, \ldots$ to be the squares of the form $[k,k+1)\times [m,m+1)$, with $k$ and $m$ integers, that are consecutively traversed by the halfline, see Figure \ref{squares}. Choosing an index $N$, there exists $M\in\mathbb{R}$ such that
$$
\bigcup_{k=1}^N S_k \cap \left\{(x,\alpha x):x\geq 0\right\} = \left\{(x,\alpha x): 0\leq x< M \right\}.
$$

Taking fractional parts in both coordinates can be seen as mapping each of the squares $S_k$ to $[0,1)^2$. Therefore, doing this for the above set gives
\begin{equation}\label{A++}
\left\{(\left\{x\right\},\left\{\alpha x\right\}):0\leq x <
M\right\} = [0,1)^2 \cap \bigcup_{k=1}^N \left\{(x,\alpha x
+y_k):x\in\mathbb{R}\right\}
\end{equation}
for numbers $y_k$ defined by the recursion
\begin{eqnarray}\label{recursion}
y_1&=&0,\nonumber\\
y_{k+1} &=& \left\{
\begin{array}{ll}
y_k+\alpha & \hspace{2cm}\textrm{if}\quad y_k < 1-\alpha,\\
y_k-1      & \hspace{2cm}\textrm{if}\quad y_k > 1-\alpha.
\end{array}
\right.
\end{eqnarray}

We will denote the set in (\ref{A++}) by $A^{++}$. The $++$ superscript reflects the fact that we started with a halfline in the first quadrant, so both coordinates are positive. Doing similar operations to halflines in the second, third and fourth quadrant,
we can define sets $A^{-+}$, $A^{--}$ and $A^{+-}$ respectively as follows:
\begin{eqnarray*}
A^{-+} &=& \left\{(1-\left\{x\right\},\left\{\alpha
x\right\}):0\leq x < M\right\}, \\
A^{--} &=& \left\{(1-\left\{x\right\},1-\left\{\alpha
x\right\}):0\leq x < M\right\}, \\
A^{+-} &=& \left\{(\left\{x\right\},1-\left\{\alpha
x\right\}):0\leq x < M\right\}.
\end{eqnarray*}

Taking the union of these four sets and intersecting with $[0,\frac{1}{2}]$ gives us a billiard orbit, as is proved in the lemma below.

\begin{lemma}\label{lemma_orbitopen}
The billiard orbit $B_{[0,M)}^\alpha$ satisfies
$$
B_{[0,M)}^\alpha = \bigcup_{u,v\in\left\{+,-\right\}} A^{uv} \cap [0,\frac{1}{2}]^2.
$$
\end{lemma}

\textbf{Proof}
Observe that
\begin{eqnarray*}
(||x||,||\alpha x||) &=& \Bigl(\min\Bigl\{\left\{x\right\},1-\left\{x\right\}\Bigr\},\min\Bigl\{\left\{\alpha x\right\},1-\left\{\alpha x\right\}\Bigr\}\Bigr)\\
&=& [0,\frac{1}{2}]^2\cap\bigcup_{a\in\left\{\left\{x\right\},1-\left\{x\right\}\right\}} \bigcup_{b\in\left\{\left\{\alpha
x\right\},1-\left\{\alpha
x\right\}\right\}} (a,b),
\end{eqnarray*}
and now take the union over all $x\in [0,M)$.\hfill $\Box$
\\
\\
In the next lemma expressions similar to (\ref{A++}) are
derived for $A^{-+}$, $A^{--}$ and $A^{+-}$.

\begin{lemma}\label{lemma_Aexpressions}
Let $y_{-k} = 1-\alpha-y_k$ for $k = 1,2,\ldots,N$. Then
\begin{eqnarray*}
A^{-+} &=& (0,1]\times[0,1) \cap \bigcup_{k=1}^N \left\{(x,-\alpha x+1-y_{-k}):x\in\mathbb{R}\right\},\\
A^{--} &=& (0,1]^2 \cap \bigcup_{k=1}^N \left\{(x,\alpha x + y_{-k}):x\in\mathbb{R}\right\}, \\
A^{+-} &=& [0,1)\times(0,1] \cap \bigcup_{k=1}^N \left\{(x,-\alpha
x+1-y_k):x\in\mathbb{R}\right\},
\end{eqnarray*}
\end{lemma}

\textbf{Proof} Define the functions
$f,g,h:\mathbb{R}^2\rightarrow\mathbb{R}^2$ by $f((x,y)) = (1-x,y)$, $g((x,y)) = (1-x,1-y)$ and $h((x,y)) = (x,1-y)$.
Applying these functions to the left hand side of (\ref{A++}), we
get $f(A^{++}) = A^{-+}$, $g(A^{++}) = A^{--}$ and $h(A^{++}) = A^{+-}$. On the other hand,
\begin{eqnarray*}
f(\left\{(x,\alpha x +y_k):x\in\mathbb{R}\right\}) &=&
\left\{(1-x,\alpha x +y_k):x\in\mathbb{R}\right\}\\
&=& \left\{(x,\alpha (1-x) +y_k):x\in\mathbb{R}\right\}\\
&=& \left\{(x,-\alpha x +1-y_{-k}):x\in\mathbb{R}\right\},
\end{eqnarray*}
whence application of $f$ to the right hand side of (\ref{A++})
leads to
\begin{eqnarray*}
f\Bigl([0,1)^2 &\cap& \bigcup_{k=1}^N \left\{(x,\alpha x
+y_k):x\in\mathbb{R}\right\}\Bigr) \\
&=&  f\Bigl([0,1)^2\Bigr) \cap
f\Bigl(\bigcup_{k=1}^N
\left\{(x,\alpha x +y_k):x\in\mathbb{R}\right\}\Bigr)\\
&=& (0,1]\times[0,1)\cap \bigcup_{k=1}^N
f\Bigl(\left\{(x,\alpha x +y_k):x\in\mathbb{R}\right\}\Bigr)\\
&=& (0,1]\times[0,1)\cap \bigcup_{k=1}^N
\left\{(x,-\alpha x +1-y_{-k}):x\in\mathbb{R}\right\},
\end{eqnarray*}
so for $A^{-+}$ we established the equality claimed in the lemma. The other two equalities for $A^{--}$ and $A^{+-}$ follow from a similar reasoning since
\begin{eqnarray*}
g(\left\{(x,\alpha x +y_k):x\in\mathbb{R}\right\}) &=&
\left\{(1-x,1-\alpha x -y_k):x\in\mathbb{R}\right\}\\
&=& \left\{(x,1-\alpha (1-x) -y_k):x\in\mathbb{R}\right\}\\
&=& \left\{(x,\alpha x +y_{-k}):x\in\mathbb{R}\right\},
\end{eqnarray*}
and
\begin{eqnarray*}
h(\left\{(x,\alpha x +y_k):x\in\mathbb{R}\right\}) &=&
\left\{(x,1-\alpha x -y_k):x\in\mathbb{R}\right\}.
\end{eqnarray*}
\hfill $\Box$

The numbers $y_k$ and $y_{-k}$ satisfy a nice relation, as is shown in the following lemma.

\begin{lemma}\label{lemma_rotation}
Let $y_0 = -\alpha$. Then the numbers $y_{-N},\ldots,y_N$ form a rotation orbit on the interval $[-\alpha,1]$. They are given by
\begin{equation}\label{y_k}
y_k = (1+\alpha)\left\{\frac{k\alpha}{1+\alpha}\right\}-\alpha  \quad \textrm{for}\quad -N\leq k\leq N.
\end{equation}
\end{lemma}

\textbf{Proof} The recursion (\ref{recursion}) can be rewritten as
$$
y_{k+1} = (y_k+2\alpha \textrm{\ mod}(1+\alpha))-\alpha,
$$
and therefore
$$
\frac{y_{k+1}+\alpha}{1+\alpha} = \frac{y_k+2\alpha}{1+\alpha}\textrm{\ mod\ }1 = \left\{\frac{y_k+2\alpha}{1+\alpha}\right\},
$$
Letting $\tilde y_k = \frac{y_k+\alpha}{1+\alpha}, k = -N,\ldots,N$ and $\tilde \alpha = \frac{\alpha}{1+\alpha}$, this reduces to
$$
\tilde y_{k+1} = \left\{\tilde y_k+\tilde \alpha\right\}.
$$
Since $y_1 = 0$, we have $\tilde y_1 = \tilde \alpha$, which leads to
$$
\tilde y_k = \left\{k\tilde \alpha\right\}\quad \textrm{for}\quad k\geq 1.
$$
On the other hand, for $k\geq 1$,
\begin{eqnarray*}
\tilde y_{-k} &=& \frac{y_{-k}+\alpha}{1+\alpha} = \frac{1-\alpha-y_k+\alpha}{1+\alpha} = \frac{1+\alpha}{1+\alpha}-\frac{y_k+\alpha}{1+\alpha}\\
&=& 1-\tilde y_k = 1- \left\{k\tilde\alpha\right\} = \left\{-k\tilde\alpha\right\},
\end{eqnarray*}
since $\tilde\alpha$ is irrational. By definition we have $\tilde y_0 = 0$, and hence
$$
\tilde y_k = \left\{k\tilde\alpha\right\} \quad \textrm{for}\quad -N\leq k\leq N.
$$
Solving for $y_k$ gives the result.
\hfill$\Box$\\
\\
In Lemma \ref{lemma_orbitopen} we already derived an expression for $B_{[0,M)}^\alpha$, but this is not so easy
to analyze directly. In the next lemma we describe $B_{[0,M]}^\alpha$ as the union of two collections of lines
intersected with $[0,\frac{1}{2}]^2$. All lines in the first collection have slope $\alpha$ and all lines in the
second collection have slope $-\alpha$.

\begin{lemma}\label{lemma_orbit}
Let $l_k^+(x) = \alpha x+y_k$ and $l_k^-(x) = -\alpha x+1-y_k$. Then
$$
B_{[0,M]}^\alpha = [0,\frac{1}{2}]^2\cap \bigcup_{u\in\left\{+,-\right\}}\bigcup_{k=-N}^N \left\{(x,l_k^u(x):x\in\mathbb{R})\right\}
$$
\end{lemma}
\textbf{Proof} This lemma will be proved by taking closures in the equation in Lemma \ref{lemma_orbitopen}.
$$
\overline{B_{[0,M)}^\alpha} = B_{[0,M]}^\alpha,
$$
since $x \mapsto (||x||,||\alpha x||)$ is a continuous function from $\mathbb{R}$ to $\mathbb{R}^2$. On the other hand,
$$
\overline{[0,\frac{1}{2}]^2\cap\bigcup_{u,v\in\left\{+,-\right\}}A^{uv}} =
[0,\frac{1}{2}]^2\cap\bigcup_{u,v\in\left\{+,-\right\}}\overline{A^{uv}},
$$
and since $A^{uv}$ is a finite collection of lines intersected by a `half open' unit square its closure is the same collection of lines but now intersected by the closed square $[0,1]^2$. Therefore,
\begin{equation}\label{++--}
\overline{A^{++}}\cup\overline{A^{--}} = [0,1]^2\cap\bigcup_{\begin{array}{l} k=-N\\k\neq 0\end{array}}^N\left\{(x,l_k^+(x)):x\in\mathbb{R}\right\}
\end{equation}
Now note that since $l_0^+(x) = \alpha x-\alpha$ we have
$$
[0,\frac{1}{2}]^2\cap\left\{(x,l_0^+(x)):x\in\mathbb{R}\right\} = \emptyset.
$$
Intersecting both sides of (\ref{++--}) with $[0,\frac{1}{2}]^2$
gives
\begin{equation}
[0,\frac{1}{2}]^2\cap\overline{A^{++}}\cup\overline{A^{--}} =
[0,\frac{1}{2}]^2\cap\bigcup_{k=-N}^N\left\{(x,l_k^+(x)):x\in\mathbb{R}\right\}
\end{equation}
Analogously it follows that
\begin{equation}
[0,\frac{1}{2}]^2\cap\overline{A^{+-}}\cup\overline{A^{-+}} =
[0,\frac{1}{2}]^2\cap\bigcup_{k=-N}^N\left\{(x,l_k^-(x)):x\in\mathbb{R}\right\}
\end{equation}
Combination of the last two equations gives the result. \hfill $\Box$ \\

\section{Sharpness of the bounds}\label{example}

In this section we present an example in which the upper bounds of Theorem \ref{theorem_13areas} are reached. This proves sharpness of the bounds.

\begin{figure}[!h]
\begin{center}
\includegraphics*[width =10cm]{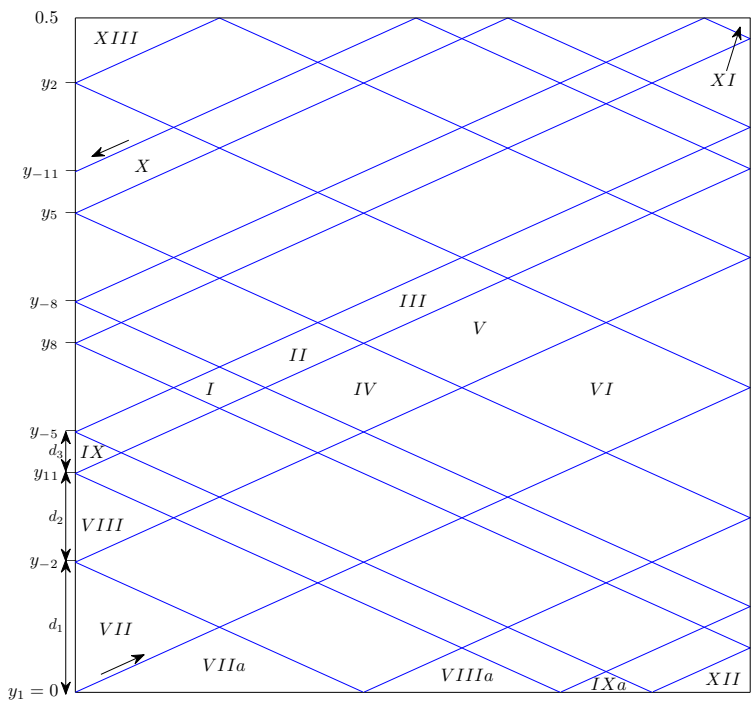}
\caption{\footnotesize Thirteen different areas, sixteen different shapes.}\label{pic_13flakes}
\end{center}
\end{figure}

\begin{example}\label{ex_13flakes}\rm{
Let $\alpha = \frac{\sqrt{10}}{7}$ and choose $N = 11$. The corresponding orbit is shown in Figure \ref{pic_13flakes}. Use Lemma \ref{lemma_rotation} to find the numbers $y_k$ and let
$$
d_1 = y_{-2}-y_1 \approx 0.0965,\quad d_2 = y_{11}-y_{-2} \approx 0.0658, \quad d_3 = y_{-5}-y_{11} \approx 0.0307
$$
denote the three different vertical distances between adjacent parallel lines. The areas of the shapes are of the following form:
\begin{equation}\label{areas}
\begin{array}{ll}
\textrm{Shapes\ $I,\ II,\ III,\ IV,\ V,\ VI$}:  &  d_id_j/2\alpha,\quad i\leq j\in\left\{1,2,3\right\}\\
\textrm{Shapes\ $VII,\ VIII,\ IX$}:          &  d_i^2/4\alpha,\quad i\in\left\{1,2,3\right\}\\
\textrm{Shape\ $X$}:                       &  d_3d_1/2\alpha - d_3^2/4\alpha\\
\textrm{Shapes\ $XI,\ XII,\ XIII$}:          &  d_i^2/8\alpha,\quad i\in\left\{1,2,3\right\}\\
\end{array}
\end{equation}
Calculating these thirteen areas indeed gives thirteen different values, where a precision of two decimals suffices.
The flakes $VII$, $VIII$ and $IX$ have the same areas as $VIIa$, $VIIIa$ and $IXa$ respectively, so the maximal number of sixteen different shapes is also reached. We checked the calculations by using the outcomes to determine the area of $[0,\frac{1}{2}]^2$}.
\end{example}

\section*{Acknowledgement}

The author thanks Michel Dekking and Cor Kraaikamp for their useful
comments.

\end{document}